\documentclass{amsproc}
\usepackage{graphicx}
\usepackage{amssymb}
\usepackage{epstopdf}
\usepackage{amsmath,amsthm,amscd,amssymb}
\usepackage{geometry}                
\numberwithin{equation}{section}

\theoremstyle{plain}
\newtheorem{theorem}{Theorem}[section]
\newtheorem{lemma}[theorem]{Lemma}

\theoremstyle{definition}
\newtheorem{definition}[theorem]{Definition}

\newtheorem{case[theorem]}{Case}

\theoremstyle{remark}

\numberwithin{equation}{section}

\def\R{\mathbb R}
\def\sd{\mathbb S^{d-1}}
\def\implies{\Longrightarrow}
\def\O{\mathcal O}

\def\zero{\mathbf 0}
\def\be{\begin{equation}}
\def\ee{\end{equation}}
\def\bes{\begin{eqnarray*}}
\def\ees{\end{eqnarray*}}
\def\st{\sqrt{3}}

\def\R{\Bbb R}

\begin{document}

\title{Equilateral triangles in subsets of ${\Bbb R}^d$ of large Hausdorff dimension}


\author{Alex Iosevich and Bochen Liu}

\date{today}

\address{Department of Mathematics, University of Rochester, Rochester, NY}

\email{iosevich@math.rochester.edu}
\email{bochen.liu@rochester.edu}

\thanks{This work was partially supported by the NSA Grant H98230-15-1-0319}

\begin{abstract} We prove that subsets of ${\Bbb R}^d$, $d \ge 4$ of large enough Hausdorff dimensions contain vertices of an equilateral triangle. It is known that additional hypotheses are needed to assure the existence of equilateral triangles in two dimensions (see \cite{CLP14}). We show that no extra conditions are needed in dimensions four and higher. The three dimensional case remains open. 

Some interesting parallels exist between the triangle problem in Euclidean space and its counter-part in vector spaces over finite fields. We shall outline these similarities in hopes of eventually achieving a comprehensive understanding of this phenomenon in the setting of locally compact abelian groups. \end{abstract}

\maketitle


\section{Introduction}

\vskip.125in

An old and classical problem that arises in different forms in geometric combinatorics, geometric measure theory, ergodic theory and other areas is to show that a sufficiently large set contains vertices of a given geometric configuration. In the realm of positive Lebesgue density, this idea can be encapsulated in the following theorem due to Tamar Ziegler, building on previous results due to Bourgain, Furstenberg, Katznelson, Weiss and others (see, for example, \cite{B86} and \cite{FKW90}). 

\begin{theorem} (T. Ziegler (\cite{Z06}) \label{zieglertheorem} Let $E \subset {\Bbb R}^d$, of positive upper Lebesgue density in the sense that 
$ \limsup_{R \to \infty} \frac{{\mathcal L}^d \{E \cap {[-R,R]}^d \}}{{(2R)}^d}>0$, where ${\mathcal L}^d$ denotes the $d$-dimensional Lebesgue measure. Let $E_{\delta}$ denote the $\delta$-neighborhood of $E$. Let $V=\{ {\bf 0}, v^1, v^2, \dots, v^{k-1}\} \subset {\Bbb R}^d$, where $k \ge 2$ is a positive integer. Then there exists $l_0>0$ such that for any $l>l_0$ and any $\delta>0$ there exists $\{x^1, \dots, x^k\} \subset E_{\delta}$ congruent to $lV=\{ {\bf 0}, lv^1, \dots, lv^{k-1}\}$. \end{theorem} 

This result nearly settles the issue of simplexes in sets of positive upper Lebesgue density, though even there an interesting open question of whether the 
$\delta$-neighborhood fudge factor in Theorem \ref{zieglertheorem} can be eliminated in the case, say, of non-degenerate triangles. If the triangle is allowed to be degenerate, an example due to Bourgain (\cite{B86}) shows that the result is not in general true without the fudge factor. 

A natural question that arises at this point is whether a compact set in ${\Bbb R}^d$ of Hausdorff dimension $s_0<d$ contains a given geometric configuration. This question is already fascinating in dimension $1$. An example due to Keleti shows that there exists a subset of $[0,1]$ of Hausdorff dimension $1$ which does not contain any arithmetic progressions of length three. However, a result due to Laba and Pramanik (\cite{LP09}) shows that there exists $s_0<1$ such that a subset of $[0,1]$ of Hausdorff dimension $s_0$ contains a progression of length three if it satisfies additional structural assumptions. 

A similar difficulty arises in higher dimensions. An example due to Falconer (\cite{Falc83}) and (independently) Maga (\cite{Ma11}) shows that there exists a set of Hausdorff dimension $2$ in ${\Bbb R}^2$ which does not contain vertices of an equilateral triangle. Once again, a result can be established with additional assumptions on the structure of the set and this was accomplished by Chung, Laba and Pramanik who proved the following, rather general result. 

\begin{definition}\label{clpConfiguration} 
Fix integers $n\geq 2$, $p\geq 3$, and $m= n\lceil \frac{p+1}{2} \rceil$.  Suppose $B_1, \dots, B_p$ are $n \times (m-n)$ matrices.

(a)  We say that $E$ contains a $p-$point $\mathcal{B}-$configuration if there exists vectors
$z\in \mathbb{R}^n $ and $w\in \mathbb{R}^{m-n}\backslash \vec{0}$ such that $\{z + B_j w \}_{j=1}^p \subset E.$ 

(b) Moreover, given any finite collection of subspaces $V_1,\dots, V_q \subset \mathbb{R}^{m-n}$ with $dim(V_i) < m-n$, we say that $E$ contains a non-trivial $p-$point $\mathcal{B}-$configuration with respect to $(V_1,\dots, V_q)$ if $\exists$ 
$z\in \mathbb{R}^n$ and $w\in \mathbb{R}^{m-n}\backslash \cup_{i=1}^{q}V_i$ such that $\{z + B_j w \}_{j=1}^p \subset E$. 

(c) Fix integers $n\geq 2$, $p\geq 3$, and $m= n\lceil \frac{p+1}{2} \rceil$. We say that a set of $n\times (m-n)$ matrices $\{ B_1, \dots, B_p\}$ is non-degenerate if  
\[rank \left( \begin{array}{c}
B_{i_2}-B_{i_1}\\
\vdots\\
B_{i_{m/n}}-  B_{i_1}\\
\end{array} \right)=m-n \] for any distinct indices $i_1,\dots,i_{ m/n} \in \{1,\dots,p\}$.
\end{definition}

\begin{theorem} \label{clp} [Chan, \L aba and Pramanik] 
Fix integers $n\geq 2$, $p\geq 3$, and $m= n\lceil \frac{p+1}{2} \rceil$.  Let $\{B_1, \dots, B_p\}$ be a collection of $n \times (m-n)$ non-degenerate matrices in the sense of part c) of Definition \ref{clpConfiguration}.  Then for any constant $C$, there exists a positive number $\epsilon_0(C,n,p,B_1,\dots,B_p) <<1$ with the following property:  Suppose the set $E \subset \mathbb{R}^n$ with $\left|E \right|=0$ supports a positive, finite, Radon measure $\mu$ with two conditions:
(a) (ball condition)  $sup_{\stackrel{x\in E}{ 0<r<1}} \frac{\mu(B(x,r)}{r^{\alpha}} \le C$ if $n-\epsilon_0 <\alpha < n$,
(b) (Fourier decay) $sup_{\xi \in \mathbb{R}^n } |\widehat{\mu}(\xi)| (1+ |\xi|)^{\beta/2} \le C.$ Then (i) $E$ contains a $p-$point $\mathcal{B}-$configuration in the sense of Definition \ref{clpConfiguration} (a). (ii) Moreover, for any finite collection of subspaces $V_1, \dots, V_q \subset \mathbb{R}^{m-n}$ with $dim(V_i) < m-n$, $E$ contains a non-trivial $p-$point $\mathcal{B}-$configuration with respect to $(V_1, \dots, V_q)$ in the sense of Definition \ref{clpConfiguration} (b). \end{theorem}

As the reader can check Theorem \ref{clp} recovers equilateral triangles in sets of Hausdorff dimension sufficiently close to $2$ under an additional structural assumption on sizes of balls and the decay rate of the Fourier transform. A natural question that arises at this point is what happens with triangles in dimensions $3$ and higher. The Chan-Laba-Pramanik result does not address this issue, even under additional structural assumptions. 

Before we state our main result, we would like to point out some interesting parallels between the configuration problems in Euclidean space and those in vector spaces over finite fields. Let ${\Bbb F}_q^d$ denote the finite with $q$ elements and let ${\Bbb F}_q^d$ denote the $d$-dimensional vector space over ${\Bbb F}_q$. Once again we may ask whether $E \subset {\Bbb F}_q^d$ contains vertices of an equilateral triangle if $E$ is sufficiently large. In two dimensions there is a fundamental arithmetic obstruction, pointed out in \cite{BIP14}, namely that if ${\Bbb F}_q$ does not contain $\sqrt{3}$, then ${\Bbb F}_q^d$ does not contain any equilateral triangles. It is not currently known whether a "metric" rather an "arithmetic" obstruction also exists. On the other hand, when $\# E \ge Cq^{\frac{2}{3}d+1}$, Hart and the first listed author proved in (\cite{HI08}) that $E$ contains vertices of every possible congruence class of triangles, including equilateral triangles with every possible side-length (in ${\Bbb F}_q$). As the reader shall see in a moment, a very similar picture emerges in Euclidean space, even though the methods and technical formulations are very different. 

\vskip.125in 

\subsection{Statement of results} The main result  of this paper is the following. We recover equilateral triangles in sets of sufficiently large Hausdorff dimension without any additional assumptions on Fourier decay rate of the underlying measure. 

\begin{theorem} \label{main} Let $E$ be a compact subset of ${\Bbb R}^d$, $d \ge 4$ and $\mu$ is a probablity Frostman measure on $E$ with $\mu(B(x,r))\leq c_\mu r^s$ for all $x\in\R^d, r>0$. Then there exists $s_0(d,c_\mu )<d$ such that if the $s>s_0$, then $E$ contains vertices of an equilateral triangle.
\end{theorem} 

\subsection{Structure of the proof} 

Our proof consists of three basic steps: \begin{itemize} 

\vskip.125in 

\item STEP 1: To construct a natural measure on the set of equilateral triangles of a given side-length with endpoints in $E$ and prove that it is finite. 

\vskip.125in 

\item STEP 2: To modify the Chang-Laba-Pramanik argument to show that with the first step as input, the integral of the measure is strictly positive.  

\vskip.125in

\item STEP 3: To show that this measure is not supported on the trivial "singleton" equilateral triangles $\{(x,x,x): x \in E \}$, thus guaranteeing the existence of a non-trivial equilateral triangle.

\end{itemize} 

\vskip.125in

{\bf{Notation.}} 
Throughout this paper, $A\lesssim B$ means there exists a constant $C_{d,c_\mu}$ such that $A\leq C_{d,c_\mu} B$.

\vskip.25in 

\section{Proof of Theorem \ref{main}} 

\subsection{Proof of STEP 1} Denote by $\sigma$ the surface measure of the surface 
$$ \{ (x,y) \in {\Bbb R}^{2d}: |x|=|y|=|x-y|=1 \}.$$

\vskip.125in

Define the measure $\nu$ on 
$$\{(x,y,z)\in E\times E\times E: |x-y|=|y-z|=|x-z|\}$$
by
\begin{equation}\label{measure}
  \begin{aligned}
   d\nu=\lim_{\delta\rightarrow 0} \mu_\delta(z)\mu_\delta(z+tx)\mu_\delta(z+ty)\,t^{d-1}\, dz\,d\sigma(x,y)\,dt,
  \end{aligned}
\end{equation}
if the limit exists, where $\mu_\delta=\mu*\phi_\delta$, $\phi_\delta=\delta^{-d}\phi(\frac{\cdot}{\delta})$ and $\phi\in C_0^\infty$ is supported in the unit ball with $\int \phi=1$. The next result shows that the limit always makes sense for a sequence tending to zero if the Hausdorff dimension of the support of $\mu$ is sufficiently large. 

\vskip.125in

\begin{theorem} \label{existence} (EXISTENCE)
If $s>\frac{2}{3}d+1$, there exists a sequence $\delta_j\rightarrow 0$ such that the limit in \eqref{measure} exists.
\end{theorem}

It suffices to show that 
$$\iiint \mu_\delta(z)\mu_\delta(z+tx)\mu_\delta(z+ty)\,t^{d-1}\, dz\,d\sigma(x,y)\,dt $$
is bounded above by a constant independent in $\delta$.

\vskip.125in

By Plancherel,

\begin{equation*}
\begin{aligned}
&\iiint \mu_\delta(z)\mu_\delta(z+tx)\mu_\delta(z+ty)\,t^{d-1}\, dz\,d\sigma(x,y)\,dt\\=&\int\cdots\int \widehat{\mu_\delta}(\xi)\widehat{\mu_\delta}(\eta)\widehat{\mu_\delta}(\zeta)e^{2\pi i((z,z,z)+(tx,ty,0))\cdot(\xi,\eta,\zeta)}\,dz\,d\sigma(x,y)\,t^{d-1}dt\,d\xi\,d\eta\,d\zeta\\=
&\int\cdots\int \widehat{\mu_\delta}(\xi)\widehat{\mu_\delta}(\eta)\widehat{\mu_\delta}(\zeta)\left(e^{2\pi i(z,z,z)\cdot(\xi,\eta,\zeta)}\,dz\right)\widehat{\sigma}(-t\xi,-t\eta)\,t^{d-1}dt\,d\xi\,d\eta\,d\zeta.
\end{aligned}
\end{equation*}

\vskip.125in

As a distribution,
$$\int e^{2\pi i(z,z,z)\cdot(\xi,\eta,\zeta)}\,dz=\delta(\xi+\eta+\zeta), $$ where $\delta$ denotes the $\delta$-distribution at the origin, hence the integral above equals

\begin{equation}\label{Fourier}
\begin{aligned}
&\iiint\widehat{\mu_\delta}(\xi)\widehat{\mu_\delta}(\eta)\widehat{\mu_\delta}(-\xi-\eta)\widehat{\sigma}(-t\xi,-t\eta)\,d\xi\,d\eta\ t^{d-1}dt\\=&\iiint\widehat{\mu}(\xi)\widehat{\phi}(\xi)\widehat{\mu}(\eta)\widehat{\phi}(\eta)\widehat{\mu}(-\xi-\eta)\widehat{\phi}(-\xi-\eta)\widehat{\sigma}(-t\xi,-t\eta)\,d\xi\,d\eta\ t^{d-1}dt.
\end{aligned}
\end{equation}

Since $\mu(B(x,r))\lesssim r^s$, we have
\begin{equation}\label{square} \int_{|\xi|\approx R} |\widehat{\mu}(\xi)|^2\,d\xi\lesssim R^{d-s}. \end{equation}

With \eqref{square}, to prove Theorem \ref{existence}, we shall show if $s>\frac{2}{3}d+1$,
\begin{equation}\label{goal}
\iiint_{|\xi|,|\eta|>R}  |\widehat{\mu}(\xi)||\widehat{\mu}(\eta)||\widehat{\mu}(\xi+\eta)||\widehat{\sigma}(t\xi,t\eta)|\,d\xi\,d\eta\ t^{d-1}dt \lesssim t^{-d+\frac{3}{2}} R^{-\frac{3s-2d-3}{2}}.
\end{equation}

\vskip.125in

The proof of \eqref{goal} requires the following stationary phase estimate. The proof will be given later.

\begin{lemma}\label{sigmahat}
Suppose $|\xi|\approx |\eta|$, then
\begin{equation}
|\widehat{\sigma}(\xi,\eta)|\lesssim
|\xi+g_{\frac{\pi}{3}}\eta|^{-\frac{1}{2}} |\xi|^{-(d-2)}\left(\sin (<\xi,\eta>)\right)^{-\frac{d-2}{2}}.
\end{equation}
where $g_{\frac{\pi}{3}}\in O(d)$ is some rotation by $\frac{\pi}{3}$ and  $<\xi,\eta>$ denotes the angle between $\xi$ and $\eta$.
\end{lemma}
\vskip.125in
Assuming $|\xi|+|\eta|\approx R2^j$, then at least two of $|\xi|, |\eta|, |\xi+\eta|$ are $\approx R 2^j$.
\vskip.125in

Case 1: $|\xi|\approx|\eta|\approx R 2^{j}$. 

\vskip.125in

By Lemma \ref{sigmahat},
\begin{equation}\label{decay1}
\begin{aligned}
&\iint_{|\xi|\approx |\eta|\approx R2^k}
  |\widehat{\mu}(\xi)||\widehat{\mu}(\eta)||\widehat{\mu}(\xi+\eta)||\widehat{\sigma}(t\xi,t\eta)|\,d\xi\,d\eta\\
  \lesssim &t^{-d+\frac{3}{2}}R^{-d+2}2^{-(d-2)j}\iint_{|\xi|\approx |\eta|\approx R2^j}
  |\widehat{\mu}(\xi)||\widehat{\mu}(\eta)||\widehat{\mu}(\xi+\eta)||\xi+g_{\frac{\pi}{3}}\eta|^{-\frac{1}{2}}\sin(<\xi,\eta>)^{-\frac{d-2}{2}}\,d\xi\,d\eta
\end{aligned}
\end{equation}

Fix $\eta$,

\begin{equation*}
\begin{aligned}
&\int_{|\xi|\approx R2^j}
|\widehat{\mu}(\xi+\eta)||\xi+g_{\frac{\pi}{3}}\eta|^{-\frac{1}{2}}\sin(<\xi,\eta>)^{-\frac{d-2}{2}}\,d\xi\\\leq &\left(\int_{|\xi|\approx R2^j}
|\widehat{\mu}(\xi+\eta)|^2\,d\xi\right)^\frac{1}{2}\left(\int_{|\xi|\approx R2^j}
  |\xi+g_{\frac{\pi}{3}}\eta|^{-1}\sin(<\xi,\eta>)^{-(d-2)}\,d\xi\right)^\frac{1}{2}\\\lesssim & (R2^j)^{\frac{d-s+\epsilon}{2}}\left(\int_{|\xi|\approx R2^j}
  |\xi+g_{\frac{\pi}{3}}\eta|^{-1}\sin(<\xi,\eta>)^{-(d-2)}\,d\xi\right)^\frac{1}{2}.
 \end{aligned}
 \end{equation*}
 
\vskip.125in

We shall need the following stationary phase argument that will be proven later on. 

\begin{lemma}\label{int}
$$\int_{|\xi|\approx R2^j}
  |\xi+g_{\frac{\pi}{3}}\eta|^{-1}\sin(<\xi,\eta>)^{-(d-2)}\,d\xi\lesssim (R2^j)^{-1+d}.$$
\end{lemma}

\vskip.125in 

Lemma \ref{int} implies that for each fixed $\eta$ with $|\eta|\approx R 2^j$,
$$\int_{|\xi|\approx R2^j}
|\widehat{\mu}(\xi+\eta)||\xi+g_{\frac{\pi}{3}}\eta|^{-\frac{1}{2}}\sin(<\xi,\eta>)^{-\frac{d-2}{2}}\,d\xi\lesssim (R2^j)^{\frac{2d-s-1}{2}}. $$

Similarly, for each fixed $\xi$ with $|\xi|\approx R 2^j$,
$$\int_{|\eta|\approx R2^j}
|\widehat{\mu}(\xi+\eta)||\xi+g_{\frac{\pi}{3}}\eta|^{-\frac{1}{2}}\sin(<\xi,\eta>)^{-\frac{d-2}{2}}\,d\eta\lesssim (R2^j)^{\frac{2d-s-1}{2}}. $$

Therefore, by Shur's test, \eqref{decay1} is bounded by 
\begin{equation}
\begin{aligned}
 &t^{-d+\frac{3}{2}}\sum_{j>k} R^{-d+2} 2^{-(d-2)j}\cdot (R2^j)^{\frac{2d-s-1}{2}} \int_{|\xi|\approx R2^j}
|\widehat{\mu}(\xi)|^2\,d\xi\\\lesssim &t^{-d+\frac{3}{2}} R^{-\frac{3s-2d-3}{2}}\sum_{j}2^{-\frac{3s-2d-3}{2}}\\\lesssim&t^{-d+\frac{3}{2}} R^{-\frac{3s-2d-3}{2}}.
\end{aligned}
\end{equation}
if $s>\frac{2}{3}d+1$.
\vskip.125in
Case 2: $|\xi|\approx |\xi+\eta|\approx R 2^j$.
\vskip.125in
Pick any pair $(x^0, y^0)$ such that $\Delta_{x^0Oy^0}$ is an equilateral triangle. Notice
\begin{equation} \label{inv}
\int f(x,y)\,d\sigma(x,y)=\int_{O(d)} f(gx^0, gy^0)\,dg,
\end{equation}
where $O(d)$ is the orthogonal group, $dg$ is the probablity Haar measure on $O(g)$ and it's independent in the choice of $x^0, y^0$. Thus
$$ \widehat{\sigma}(\xi,\eta)= \int e^{-2\pi i (gx^0\cdot\xi+gy^0\cdot\eta)}\,dg= \int e^{-2\pi i (gx^0\cdot\xi+gy^0\cdot\eta)}\,dg=\int e^{-2\pi i (g(y^0-x^0)\cdot(-\xi)+gy^0\cdot(\xi+\eta))}\,dg.$$
Since $\Delta_{x^0Oy^0}$ is an equilateral triangle implies that $\Delta_{y^0O(y^0-x^0)}$ is also an equilateral triangle,
$$ \widehat{\sigma}(\xi,\eta)=\widehat{\sigma}(-\xi,\xi+\eta).$$

With $\zeta=\xi+\eta$
\begin{equation}
\begin{aligned}
&\iint_{|\xi|\approx |\xi+\eta|\approx R2^k}
  |\widehat{\mu}(\xi)||\widehat{\mu}(\eta)||\widehat{\mu}(\xi+\eta)||\widehat{\sigma}(t\xi,t\eta)|\,d\xi\,d\eta\\=&\iint_{|\xi|\approx |\zeta|\approx R2^k}
  |\widehat{\mu}(\xi)||\widehat{\mu}(\zeta-\xi)||\widehat{\mu}(\zeta)||\widehat{\sigma}(-t\xi,t\zeta)|\,d\xi\,d\zeta\\\lesssim &t^{-d+\frac{3}{2}} R^{-\frac{3s-2d-3-3\epsilon}{2}},
\end{aligned}
\end{equation}
by the argument in Case 1.

\vskip.125in 

The discussion in case 1,2 implies \eqref{goal} and completes the proof of Theorem \ref{existence}.
\vskip.25in

\subsection{STEP 2}
In \cite{CLP14}, the authors proved the following theorem
\begin{theorem}[Proposation 5.1, \cite{CLP14}]\label{CLP14}
Let

$$\Lambda(f)=  \int_{\R^n}\int_{\R^{n-m}} \prod_{j=1}^k f(z+B_jx)\,dz\,dx, $$
where $\{B_1,\dots,B_k\}$ is non-degenerate (see Definition \ref{clpConfiguration} above). Then for every $\lambda, M>0$, there exists a constant $c(\lambda, M)>0$ with the following property: for every function $f:[0,1]^n\rightarrow \R$, $0\leq f\leq M$, $\int f\geq\lambda$, we have $\Lambda(f)\geq c(\lambda,M)$.
\end{theorem}
\vskip.125in
It's not hard to check that
$$\iiint \mu_\delta(z)\mu_\delta(z+tx)\mu_\delta(z+ty)\,t^{d-1}\, dz\,d\sigma(x,y)\,dt =\int\int_{\R^d}\int_{\R^d} \mu_\delta(z)\mu_\delta(z+x)\mu_\delta(z+gx)\,dz\,dx\,dg,$$
where $dg$ is the measure on the set of $\frac{\pi}{3}$-rotations induced by the map $$e_1\rightarrow \{x\in\R^d:\Delta_{xOe_1}\text{\ is an equilateral triangle}\}=S^{d-2}.$$
\vskip.125in
Note that $||\mu_\delta||_\infty\lesssim \delta^{s-d}$, $\int \mu_\delta=1$.
\vskip.125in
Fix $g$, applying Theorem \ref{CLP14} with $n=2m=2d, k=3$, $B_1=0, B_2=id, B_3=g$, it follows that
$$\iint \mu_\delta(z)\mu_\delta(z+x)\mu_\delta(z+gx)\,dz\,dx\gtrsim c(\delta^{s-d}). $$

Since the set of $\frac{\pi}{3}$-rotations is compact, after integrating in $g$ we have
$$I:=\iiint \mu_\delta(z)\mu_\delta(z+x)\mu_\delta(z+gx)\,dz\,dx\,dg\gtrsim c(\delta^{s-d}). $$

The proof in section 2.1 implies
$$ \int d\nu= \iiint  \widehat{\mu}(\xi)\widehat{\mu}(\eta)\widehat{\mu}(-\xi-\eta)\widehat{\sigma}(-t\xi,-t\eta)\,d\xi\,d\eta\ t^{d-1}dt.$$

Therefore, by \eqref{Fourier}, \eqref{goal},
$$\left|\int d\nu-I\right|\lesssim \iiint_{|\xi|,|\eta|>\delta^{-1}}  |\widehat{\mu}(\xi)||\widehat{\mu}(\eta)||\widehat{\mu}(\xi+\eta)||\widehat{\sigma}(t\xi,t\eta)|\,d\xi\,d\eta\ t^{d-1}dt \lesssim \delta^{\frac{3s-2d-3}{2}}.$$
\vskip.125in
Let $\delta=e^{\frac{1}{s-d}}$. It follws that $I\gtrsim 1$ while 
$$\left|\int dv-I\right|\rightarrow 0$$
as $s\rightarrow d$, which implies that there exists $s_0(d,c_\mu)<d$ such that if $s>s_0(d,c_\mu)$, 
$$\int d\nu>0.$$
\vskip.25in
\subsection{Step 3}
The estimate \eqref{goal} shows that if $s>\frac{2}{3}d+1$, as a function of $t$,
$$\iint |\widehat{\mu}(\xi)||\widehat{\mu}(\eta)||\widehat{\mu}(\xi+\eta)||\widehat{\sigma}(t\xi,t\eta)|\,d\xi\,d\eta\ t^{d-1}\in L^1([0,1]).$$

Dominated convergence theorem and the discussion in Step 2 implies that if $s>s_0(d,c_\mu)$, there exists $t_0>0$ and a sequence $\delta_j\rightarrow 0$ such that
$$\lim_{\delta_j\rightarrow 0} \int_{t_0}^1\iint\mu_\delta(x)\mu_\delta(x+ty)\mu_\delta(x+tz)\,t^{d-1}\, dx\,d\sigma(y,z)\,dt$$
is well defined and positive. This means there exists a positive measure
$$ \lim_{\delta_j\rightarrow 0} \left(\int_{t_0}^1\mu_\delta(x)\mu_\delta(x+ty)\mu_\delta(x+tz)\,t^{d-1}\,dt\right)\, dx\,d\sigma(y,z) $$
supported on 
$$\{(x,y,z)\in E\times E\times E: |x-y|=|y-z|=|x-z|>t_0\}, $$
 which completes the proof of Theorem \ref{main}.\vskip.25in

\section{Proof of Lemma \ref{sigmahat}}
The upper bound of $|\widehat{\sigma}(\xi,\eta)|$ was first studied in \cite{GGIP15} (Lemma 3.19), where the authors obtain 
$$ |\widehat{\sigma}(\xi,\eta)|\lesssim (1+|\xi|+|\eta|)^{-\frac{d-1}{2}}.$$
\vskip.125in
We follow the idea of the proof in the \cite{GGIP15}, with more delicate computation and observation. 
\vskip.125in
By partition of unity and the rotation invariance of $\sigma$ (see \eqref{inv}), one can only consider a small neighborhood of $x^0=(0,\zero',1)$ and $y^0=(\frac{\sqrt3}2,\zero',\frac12)$, where we write $\R^d\owns x=(x_1,x',x_d)$. Introduce local coordinates on $\sd$ near $x^0,y^0$, resp.,
\bes
x(u)&=&\left(u,1-\frac{|u|^2}2 \right)+\O(|u|^3),\, u=(u_1,u')\in\R^{d-1},\, |u|< \epsilon,\hbox{ and }\\
y(v)&=& \left(\frac{\sqrt3}2 +v_1,v',\frac12-\st v_1-4v_1^2-|v'|^2 \right)+\O(|v|^3),\, v=(v_1,v')\in\R^{d-1},\, |v|<\epsilon.
\ees
All calculations that follow will be modulo $\O^3:=\O(|u,v|^3)$. 
\vskip.125in
Notice $|x|=|y|=|x-y|=1$ implies $x\cdot y=\frac{1}{2}$. However, 
\bes
x(u)\cdot y(v)-\frac12&=&\frac{\st}2 u_1-\st v_1+u\cdot v -\frac{|u|^2}4-4v_1^2-|v'|^2 +\O^3=0\Leftrightarrow \\
(\st-u_1)v_1 + 4v_1^2&=&\frac{\st}2 u_1+u'\cdot v'-\frac{|u|^2}4-|v'|^2+\O^3.
\ees
We use the  quadratic terms in the implicit function theorem in one variable,
\be
a_1 s+ a_2 s^2=t\implies s=a_1^{-1}t-a_1^{-3}a_2t^2+\O(t^3),\, s,t\searrow 0,
\ee
to solve for $v_1$ in terms of $u_1$, with $u',v'$ as parameters:
$$v_1=\frac{1}{2} u_1- \frac{\sqrt{3}}{4}u_1^2-\frac{\sqrt{3}}{12}|u'|^2-\frac{1}{\sqrt{3}}|v'|^2+\frac{1}{\sqrt{3}}u'\cdot v'+O^3.$$

Since $|x|=|y|=|x-y|=1$ is equivalent to $x,y\in S^{d-1}, x\cdot y=\frac{1}{2}$, a neighborhood of $(x^0,y^0)$ can be parametrized by  
$$ \left(u_1, u', 1-\frac{|u|^2}{2}; \frac{\sqrt{3}}{2}+\frac{u_1}{2}- \frac{\sqrt{3}}{4}u_1^2-\frac{\sqrt{3}}{12}|u'|^2-\frac{1}{\sqrt{3}}|v'|^2+\frac{1}{\sqrt{3}}u'\cdot v',v',\frac{1}{2}-\frac{\sqrt{3}}{2}u_1-\frac{1}{4}u_1^2+\frac{1}{4}|u'|^2-u'\cdot v' \right)$$
$$:= \vec{U},$$
modulo $\O^3$.
\vskip.125in
Then the Fourier transform of the measure $\sigma$ can be written as
$$\int e^{-2\pi i \,\vec{U}\cdot (\xi, \eta)}\,du_1\,du'\,dv',$$
and, by the invariance, we may asuume $(u_1,u',v')=\zero$ is a critical point, if exists.

The gradient of the phase function $\vec{U}\cdot (\xi,\eta)$ is,
$$\frac{\partial}{\partial u_1}=\xi_1-u_1\xi_d+\frac{1}{2}\eta_1-\frac{\sqrt{3}}{2}u_1\eta_1-\frac{\sqrt{3}}{2}\eta_d-\frac{1}{2}u_1\eta_d;$$
$$\frac{\partial}{\partial u'}=\xi'-u'\xi_d-\frac{1}{2\sqrt{3}}u'\eta_1+\frac{1}{\sqrt{3}}v'\eta_1+\frac{1}{2}u'\eta_d-v'\eta_d;$$
$$\frac{\partial}{\partial v'}=-\frac{2}{\sqrt{3}}v'\eta_1+\frac{1}{\sqrt{3}}u'\eta_1+\eta'-u'\eta_d.$$

And the Hessian is 
\begin{equation*}
(-\xi_d-\frac{\sqrt{3}}{2}\eta_1-\frac{1}{2}\eta_d)\bigoplus_2^{d-1}
\begin{pmatrix}
-\xi_d-\frac{1}{2\sqrt{3}}\eta_1+\frac{1}{2}\eta_d & \frac{1}{\sqrt{3}}\eta_1-\eta_d\\ \frac{1}{\sqrt{3}}\eta_1-\eta_d & -\frac{2}{\sqrt{3}}\eta_1
\end{pmatrix}.
\end{equation*}
\vskip.125in
Since $(u_1,u',v')=\zero$ is a critical point, if exists, it follows that 
\begin{equation}\label{critical}
\xi_1+\frac{1}{2}\eta_1-\frac{\sqrt{3}}{2}\eta_d=0,\,\xi'=\eta'=0.
\end{equation}

By \eqref{critical}, $-\xi_d-\frac{1}{2\sqrt{3}}\eta_1+\frac{1}{2}\eta_d=-\xi_d+\frac{1}{\sqrt{3}}\xi_1$, $\frac{1}{\sqrt{3}}\eta_1-\eta_d=-\frac{2}{\sqrt{3}}\xi_1$. Therefore the determinant of the Hessian equals
$$\left|-\xi_d-\frac{\sqrt{3}}{2}\eta_1-\frac{1}{2}\eta_d\right|\left|\frac{2}{\sqrt{3}}\eta_1(\xi_d-\frac{1}{\sqrt{3}}\xi_1)-\frac{4}{3}\xi_1^2\right|^{d-2}$$
$$=\left|-\xi_d-\frac{\sqrt{3}}{2}\eta_1-\frac{1}{2}\eta_d\right|\left|\frac{2}{\sqrt{3}}\eta_1\xi_d-\frac{2}{3}\eta_1\xi_1-\frac{4}{3}\xi_1^2\right|^{d-2}$$
$$=\left|-\xi_d-\frac{\sqrt{3}}{2}\eta_1-\frac{1}{2}\eta_d\right|\left|\frac{2}{\sqrt{3}}\eta_1\xi_d-\xi_1\left(\frac{2}{3}\eta_1+\frac{4}{3}\xi_1\right)\right|^{d-2},$$
agian, by \eqref{critical},
$$=\left|-\xi_d-\frac{\sqrt{3}}{2}\eta_1-\frac{1}{2}\eta_d\right|\left|\frac{2}{\sqrt{3}}\eta_1\xi_d-\frac{2}{\sqrt{3}}\xi_1\eta_d\right|^{d-2},$$
since $\xi'=\eta'=0$, we are working on $\R^2$ with $\xi,\eta$, so
$$=C\left|-\xi_d-\frac{\sqrt{3}}{2}\eta_1-\frac{1}{2}\eta_d\right||\xi|^{d-2}|\eta|^{d-2}\left(\sin (<\xi,\eta>)\right)^{d-2}.$$

Putting \eqref{critical} and the first factor of the determinant of the Hessian together as a vector in $\R^2$,
\begin{equation*}
\begin{pmatrix}
\xi_1+\frac{1}{2}\eta_1-\frac{\sqrt{3}}{2}\eta_d\\
\xi_d+\frac{\sqrt{3}}{2}\eta_1+\frac{1}{2}\eta_d
\end{pmatrix}
=\xi+g_{\frac{\pi}{3}}\eta,
\end{equation*}
where $g_{\frac{\pi}{3}}\eta$ is the vector obtained by rotating $\eta$ counterclockwise by $\frac{\pi}{3}$.

This means, at the critical point, $|-\xi_d-\frac{\sqrt{3}}{2}\eta_1-\frac{1}{2}\eta_d|=|\xi+g_{\frac{\pi}{3}}\eta|$. By stationary phase,

$$|\hat{\sigma}(\xi,\eta)|\lesssim |\xi+g_{\frac{\pi}{3}}\eta|^{-\frac{1}{2}} |\xi|^{-\frac{d-2}{2}}|\eta|^{-\frac{d-2}{2}}\left(\sin (<\xi,\eta>)\right)^{-\frac{d-2}{2}},$$

where $g_{\frac{\pi}{3}}\in O(n)$ is some rotation by $\frac{\pi}{3}$.

\vskip.25in
\section{Proof of Lemma \ref{int}}

If $<\xi,\eta>\gtrsim1$ and $|\xi+g_{\frac{\pi}{3}}\eta|\gtrsim |\xi|\approx R2^j$, the lemma is trivial.
\vskip.125in
If $<\xi,\eta>$ is small, $$|\xi+g_{\frac{\pi}{3}}\eta|^{-1}\approx |\xi|^{-1}\approx (R2^j)^{-1}.$$
Using polar coordinate we see that 
$$\int_{\substack{|\xi|\approx R2^j\\ <\xi,\eta>\,\text{small}}}|\xi+g_{\frac{\pi}{3}}\eta|^{-1}
  \sin(<\xi,\eta>)^{-(d-2)}\,d\xi\approx (R2^j)^{-1+d} \int_{S^{d-1}}
  \sin(<\omega,\eta>)^{-(d-2)}\,d\omega_{d-1}=2\pi \cdot (R2^j)^{-1+d},$$
where the last equality comes from the formula $d\omega_{d-1}=(\sin\theta)^{d-2}\,d\theta\,d\omega_{d-2}$.
\vskip.125in
If $|\xi+g_{\frac{\pi}{3}}\eta|$ is small, it follows that $<\xi,\eta>\gtrsim 1$ since $|\xi|\approx|\eta|$. Thus with $\zeta=\xi+g_{\frac{\pi}{3}}\eta$
$$ \int_{\substack{|\xi|\approx R2^j\\ |\xi+g_{\frac{\pi}{3}}\eta|\,\text{small}}}|\xi+g_{\frac{\pi}{3}}\eta|^{-1}
  \sin(<\xi,\eta>)^{-(d-2)}\,d\xi\lesssim \int_{|\zeta|\lesssim R2^j}|\zeta|^{-1}\,d\zeta\lesssim (R2^j)^{-1+d}.$$

\end{document}